\documentstyle[12pt]{article}
\begin{document}

{\LARGE \bf Hausdorff Continuous Solutions of } \\

{\LARGE \bf Arbitrary Continuous Nonlinear PDEs } \\
{\LARGE \bf through the Order Completion Method } \\ \\

{\bf Elem\'{e}r E ~Rosinger} \\ \\
{\small \it Department of Mathematics \\ and Applied Mathematics} \\
{\small \it University of Pretoria} \\
{\small \it Pretoria} \\
{\small \it 0002 South Africa} \\
{\small \it eerosinger@hotmail.com} \\ \\

{\bf Abstract} \\

In 1994 we showed that very large classes of systems of nonlinear PDEs have solutions which
can be assimilated with usual measurable functions on the Euclidean domains of definition of
the respective equations. Recently, the {\it regularity} of such solutions has significantly
been improved by showing that they can in fact be assimilated with Hausdorff continuous
functions. The method of solution of PDEs is based on the Dedekind order completion of spaces
of smooth functions which are defined on the domains of the given equations. In this way, the
method does {\it not} use functional analytic approaches, or any of the customary
distributions, hyperfunctions, or other generalized functions. \\ \\

{\bf Type independent existence and regularity results for large classes of systems of
nonlinear PDEs} \\

Ten years ago, in [4], the following significant {\it threefold} breakthrough was obtained
with respect to solving large classes of nonlinear PDEs, see MR 95k:35002. Namely : \\

a) arbitrary nonlinear PDEs of the form \\

(1)~~~~~~ $ F ( x, U(x), ~.~.~.~ , D^p U(x), ~.~.~.~ ) ~=~ f(x),~~ x \in \Omega $ \\

with $F$ jointly continuous in all it arguments, $f$ in a class of measurable functions,
$\Omega \subseteq {\bf R}^n$ arbitrary open, $p \in {\bf N}^n$, with $| p |\leq m$, for $m \in
{\bf N}$ arbitrary given, and the unknown function $U : \Omega \longrightarrow {\bf R}$, were
proven to have \\

b) solutions $U$ which can be assimilated with usual measurable functions on $\Omega$, and \\

c) the solution method was based on the Dedekind order completion of suitable spaces of smooth
functions on $\Omega$. \\

In fact, the conditions at a)~ can further be relaxed by assuming that $F$ may admit certain
{\it discontinuities}, namely, that it is continuous only on $(\Omega \setminus \Sigma )
\times {\bf R}^{m^*}$, where $\Sigma$ is a closed, nowhere dense subset of $\Omega$, while
$m^*$ is the number of arguments in $F$ minus $n$. This relaxation on the continuity of $F$
may be significant since such subsets of discontinuity $\Sigma$ can have arbitrary large
positive Lebesgue measure. \\

The method of order completion and the results on the existence and regularity of solutions
can easily be extended to {\it systems} of nonlinear PDEs of the above form (1). Furthermore,
initial and/or boundary value problems can be dealt with easily by this order completion
method. \\
In this way, the solutions of the unprecedented large class of nonlinear PDEs in (1) can be
obtained {\it without} the use of any sort of distribution, hyperfunctions, generalized
functions, or of methods of functional analysis. Moreover, one obtains a general, {\it blanket
regularity}, given by the fact that the solutions constructed can be assimilated with usual
measurable functions on the corresponding domains $\Omega$ in Euclidean spaces. \\

Recently, in collaboration with R. Anguelov, see [1], a further significant improvement of the
above mentioned 1994 results was obtained. Namely, this time we can further improve the {\it
regularity} properties of the solutions by proving that they always belong to the
significantly smaller class of Hausdorff continuous functions on the open domains $\Omega$,
see Appendix for a short account on Hausdorff continuous functions. \\

It should be noted that the results in [4] on existence of solutions do for the first time in
the literature manage fully to overcome the celebrated 1957 Hans Lewy impossibility, see [5],
and in fact do so with a large nonlinear margin. \\

Also, the existence results in [4], and thus their mentioned recent improvement with resepct
to the regularity of solutions, when solving large classes of nonlinear PDEs, {\it supersede}
to a good extent the earlier similar ones obtained through the algebraic nonlinear theory of
generalized functions introduced by the author in the 1960s, and developed since then alone or
in collaborations, see 46F30 in the AMS Subject Classification
at www.ams.org/index/msc/46Fxx.html, as well as [6], [7], [3, p. 7] and the literature cited
there, or MR 89g:35001, MR 92d:46098, Zbl.Math.717 35001, Bull.AMS, Jan. 1989, Vol. 20, No. 1,
96-101. \\

To further facilitate the understanding of the above mentioned results, it may be useful to
point to the following. In his latest 2004 edition of his Springer Universitext book "Lectures
on PDEs", see [2], V I Arnold starts on page 1 with the statement : \\

"In contrast to ordinary differential equations, there is {\it no unified theory} of partial
differential equations. Some equations have their own theories, while others have no theory at
all. The reason for this complexity is a more complicated geometry ..." (italics added) \\

However, as the above mentioned results show, since [4], there is an {\it existence and
regularity} of solutions theory for the large class of systems of nonlinear PDEs of the form
in (1). Moreover, recently, the regularity result has been significantly improved by proving
the existence of Hausdorff continuous solutions for such general nonlinear systems of
PDEs. \\ \\

{\bf Appendix. Hausdorff continuous functions} \\

We shall deal with functions whose values can be usual or {\it extended} real numbers, that is,
elements in \\

$ \overline{{\bf R}} = {\bf R} \cup \{ -\infty, +\infty \} $ \\

Moreover, we shall allow the values of the functions to be not only numbers in
$\overline{{\bf R}}$, but also {\it closed intervals} of such numbers, namely \\

$ [ a, b ] \subseteq \overline{{\bf R}},~~~ a, b \in \overline{{\bf R}},~ a \leq b $ \\

It turns out to be quite surprising how much more appropriate such a framework is when one
deals with large classes of {\it nonsmooth} functions in what is usually called Real
Analysis. \\

Indeed, by considering such {\it interval valued} functions one obtains a systematic and
effective way to study and deal with a large variety of nonsmooth functions. Furthermore, one
can gain important insights into the properties of such nonsmooth functions, properties which
in fact are not available in the usual approach. \\
It appears therefore that with the emergence in the second half of the 19-th century of a
rigorous approach to Analysis, and specifically, with the Dirichlet definition of a function
as having values given only and only by one single number, a certain undesired limitation was
imposed in an unintended manner, especially what the study of nonsmooth functions is
concerned. \\
Somewhat later, towards the end of the 19-th century, Baire brought in the concepts of {\it
lower} and {\it upper semi-continuous} functions, when dealing with nonsmooth real valued
functions. And in effect, he associated with each real valued function $f$, {\it two} other
real, or extended real valued functions $I ( f )$ and $S ( f )$, with $I ( f ) \leq f \leq
S ( f )$, which proved to be particularly helpful. However, following the prevailing mentality,
each of these three functions were considered as being single valued. \\

As it turns out, however, by considering {\it interval valued} functions, such as for instance
$F ( f ) = [ I ( f ), S ( f ) ]$, one can significantly improve on the understanding and
handling of nonsmooth functions. \\
The study of interval valued functions can, among others, show that the particular case of
functions which have values given by one single number is appropriate for continuous functions
only. On the other hand, nonsmooth functions are much better described by suitably associated
interval valued functions. \\

Indeed, in the case of functions $f$ which are {\it not} continuous, a much better description
can be obtained by considering them given by a {\it pair} of usual point valued functions,
namely $f = [~ \underline{f}, \overline{f} ~]$, thus leading to interval valued functions. And
then, a natural class which replaces, and also extends, the usual point valued continuous
functions is that of {\it Hausdorff-continuous} interval valued functions. The distinctive and
{\it essential} feature of these Hausdorff-continuous functions $f = [~ \underline{f},
\overline{f} ~]$ is a condition of {\it minimality} with respect to the {\it gap} between
$\underline{f}$ and $\overline{f}$, with the further requirement that $\underline{f}$ be lower
semi-continuous, and $\overline{f}$ be upper semi-continuous. \\
In retrospect, it is surprising to see how near Baire came to such a treatment of nonsmooth
functions, what deep results he obtained, and its correspondent for lower semi-continuous
functions, and yet followed the prevailing trend which considered functions as having to have
point, and not interval values. \\

A good measure of the {\it naturaleness} of interval valued functions can be seen in the
results related to the Dedekind order completion of various spaces of continuous functions.
And it is precisely such recently obtained results which allow for the mentioned significantly
increased regularity properties of solutions of PDEs. \\

These Dedekind order completions prove to be subspaces of Hausdorff-conti-\\nuous, thus
interval valued functions. By the way, the space of Hausdorff-continuous functions itself is
order complete. \\
Such results extend easily to functions defined on large classes of topological spaces. \\
A further indication of the natural role interval values play in the study of nonsmooth
functions can be found in the Differential and Integral Calculus being presently developed for
functions with such values. \\

It will be useful to start by introducing a few notations. Let \\

(A.1) \quad $ \overline{{\bf I R}} ~=~ \{~ [\underline{a},~ \overline{a}] ~~|~~ \underline{a},~ \overline{a} \in
           \overline{{\bf R}} ~=~ {\bf R} \cup \{ - \infty, + \infty \},~~ \underline{a}~ \leq ~\overline{a}~ \} $ \\

be the set of all finite or infinite closed intervals.

The functions which we consider can be defined on arbitrary topological spaces $\Omega$. For
the purposes of the nonlinear PDEs studied in this book, however, it will be sufficient to
assume that $\Omega \subseteq {\bf R}^n$ are arbitrary open subsets. \\
Let us now consider the set of interval valued functions \\

(A2) \quad $ {\bf A} ( \Omega ) ~=~ \{~ f : \Omega ~\longrightarrow~
                                                \overline{{\bf I R}} ~\} $ \\

By identifying the point $a \in \overline{{\bf R}}$ with the degenerate interval $[a,a] \in
\overline {{\bf I R}}$, we consider $\overline{{\bf R}}$ as a subset of
$\overline{{\bf I R}}$. In this way ${\bf A} ( \Omega )$ will contain the set of functions
with extended real values, namely \\

(A3) \quad $ {\cal A} ( \Omega ) ~=~  \{~ f : \Omega ~\longrightarrow~ \overline{{\bf R}} ~\} ~\subseteq~ {\bf A}
                                     ( \Omega ) $ \\

We define a partial order $\leq$ on $\overline{{\bf I R}}$ by \\

(A4)~~~ $ [ \underline{a},~ \overline{a} ] \leq [ \underline{b},~\overline {b} ] ~~~\Longleftrightarrow~~~ \underline{a}~ \leq
~\underline{b},~~ \overline{a}~ \leq ~\overline{b} $ \\

Now on ${\bf A} ( \Omega )$ we define the partial order induced by (A2.1.4) in the usual point-wise way, namely, for
$f,~ g \in {\bf A} ( \Omega )$, we have \\

(A5)~~~ $ f \leq g ~~~\Longleftrightarrow~~~ f(x) \leq g(x),~~ x \in \Omega $ \\

Clearly, when restricted to ${\cal A} ( \Omega )$, the above partial order on ${\bf A}
( \Omega )$ reduces to the usual one among point valued functions. \\

Given an interval $a = [ \underline{a},~ \overline{a} ] \in \overline{{\bf I R}}$, we denote\\

(A6) \quad $ w(a) ~=~~  \begin{array}{| l}
                                          ~~ \overline{a} - \underline{a}~~~ \mbox{if}~~ \underline{a},~ \overline{a}~~ \mbox{finite} \\ \\
                                          ~~ \infty~~~ \mbox{if}~~~ \overline{a}= \infty~ \mbox{and}~ \underline{a}~ \mbox{finite,~ or} \\
                                                 ~~~~~~~~~~~~~\underline{a} = - \infty~ \mbox{and}~ \overline{a}~ \mbox{finite,~ or} \\
                                                 ~~~~~~~~~~~~~\underline{a} = - \infty~ \mbox{and}~ \overline{a} = \infty \\ \\
                                          ~~ 0~~~ \mbox{if}~~ \underline{a} ~=~ \overline{a} ~=~ \pm\infty
                                 \end{array} $ \\

which is called the {\it width} of the interval $a$. Also, we denote by \\

(A7) \quad $ | a | = \mbox{max} \{~ | \underline{a} |,~ | \overline{a} | ~\} $ \\

the {\it modulus} of the interval $a = [ \underline{a},~ \overline{a} ] \in \overline{{\bf I R}}$. \\

In this way \\

(A8) \quad $ {\cal A} ( \Omega ) ~=~ \{~ f \in {\bf A} ( \Omega ) ~|~ w ( f ( x ) ) = 0,~~ x \in \Omega ~\}
                              ~\subseteq~ {\bf A} ( \Omega ) $ \\

Let $f \in {\bf A} ( \Omega )$. For every $x\in\Omega$, the value of $f$ is an interval,
namely \\

$ f ( x ) ~=~ [~ \underline{f}(x),~\overline{f}(x) ~],~~~ \mbox{with}~~ \underline{f}(x),~\overline{f}(x) \in \overline{{\bf R}},~
                            \underline{f}(x) \leq \overline{f}(x) $ \\

Hence, every function $f \in {\bf A} ( \Omega )$ can be written in the form \\

(A9) \quad $ f ~=~ [~ \underline{f},~ \overline{f} ~],~~ \mbox{with}~~ \underline {f},~ \overline{f} \in
                                 {\cal A} ( \Omega ),~~\underline{f} ~\leq~ f ~\leq~ \overline{f} $ \\

and \\

(A10) \quad $ f \in {\cal A} ( \Omega ) ~~~\Longleftrightarrow~~~  \underline{f} ~=~ f ~=~ \overline{f} $ \\

In the particular case of functions in ${\cal A} ( \Omega )$, that is, with extended real, but
point, and not nondegenerate interval values, a number of results in the sequel were obtained
by Baire. \\
Most of the more general results concerning functions in ${\bf A} ( \Omega )$, that is, with
values finite or infinite closed intervals, have recently been developed by Anguelov. \\

For $x \in \Omega$, we denote by ${\cal V}_x$ the set of all neighbourhoods $V \subseteq
\Omega$ of $x$. \\
Let us consider the pair of mappings $I,~ S : {\bf A}( \Omega ) \rightarrow {\cal A}(\Omega)$, called {\it lower} and {\it
upper Baire operators}, respectively, where for every function $f \in {\bf A} ( \Omega )$, we define \\

(A11)~~~ $ I ( f ) ( x ) ~=~ \sup_{V \in {\cal V}_x}~ \inf~ \{~ z \in f(y) ~~|~~ y \in V ~\} $ \\

(A12)~~~ $ S ( f ) ( x ) ~=~ \inf_{V \in {\cal V}_x}~ \sup~ \{~ z\in f(y) ~~|~~ y \in V ~\} $ \\

In Baire, these two operators were considered and studied in the particular case of functions
$f \in {\cal A} ( \Omega )$. \\

In view of the main interest here in this book in {\it interval valued} functions $f \in {\bf A} ( \Omega )$, it is useful to
consider as well the following third mapping, namely, $F : {\bf A} (\Omega) \rightarrow {\bf A} ( \Omega )$, defined by \\

(A13)~~~ $ F ( f ) ( x ) ~=~ [~ I ( f ) ( x ),~ S ( f ) ( x ) ~],~~ f \in {\bf A} ( \Omega ),~ x \in \Omega,~ $ \\

and called the {\it graph completion operator}. \\

The lower and upper Baire operators $I$ and $S$, and consequently, the graph completion
operator $F$, applied to any interval valued function $f = [~ \underline{f},~\overline{f} ~]
\in {\bf A} ( \Omega )$ can  now be conveniently represented in terms of the functions
$\underline{f}$ and $\overline{f}$. Indeed, from (A2.1.11), (A2.1.12) it is easy to see that \\

(A14)~~~ $ I ( f ) ~=~ I ( \underline{f} ),~~~ S( f ) ~=~ S( \overline{f} ) $ \\

Hence $F( f )$ can be written in the form \\

(A15)~~~ $ F ( f ) ~=~ [~ I ( f ),~ S ( f ) ~] ~=~ [~ I ( \underline{f} ),~ S ( \overline{f} ) ~]$ \\

Let us note that for every function $f \in {\bf A} ( \Omega )$ we have the relations, see
(A5), (A9) \\

(A16)~~~ $ I ( f ) ~=~ I ( \underline{f} ) ~~\leq~~ \underline{f} ~~\leq~~ f ~~\leq~~ \overline{f} ~~\leq~~
                                          S ( \overline{f} ) ~=~ S ( f )  $ \\

and, thus, the inclusions \\

(A17)~~~ $ f ( x ) \subseteq F ( f ) ( x ),~~ x \in \Omega $ \\

Furthermore, the lower Baire operator $I : f \rightarrow I ( f )$, the upper Baire operator
$S : f \rightarrow S( f )$ and the graph completion operator $F : f \rightarrow F ( f ) =
[~ I( f ),~ S ( f ) ~]$ are all monotone with respect to the order $\leq$ in (A5) on ${\bf A}
( \Omega )$, which means that for every two functions $f, g \in {\bf A} ( \Omega )$ we have \\

(A18)~~~ $ f \leq g ~~~\Longrightarrow~~~ I ( f ) \leq I ( g ),~~ S ( f ) \leq S ( g ),~~ F ( f ) \leq F ( g ) $ \\

The operator $F$ is also monotone with respect to inclusion, namely \\

(A19)~~~ $ f ( x ) \subseteq g ( x ),~~ x \in \Omega\ ~~~\Longrightarrow~~~ F ( f ) ( x )
\subseteq F ( g ) ( x ),~~ x\in \Omega $ \\

With an immediate extension of Baire, one can also show that all three operators are
idempotent, that is, for every $f \in {\bf A} ( \Omega )$, we have \\

(A20)~~~ $ I ( I ( f ) ) ~=~ I ( f ),~~~ S ( S ( f ) ) ~=~ S ( f ),~~~
                                                 F ( F ( f ) ) ~=~ F ( f ) $ \\ \\

{\bf Definition A1} \\

A function $f \in {\bf A} ( \Omega )$ is called {\it segment-continuous}, or in short, {\it
s-continuous}, if and only if \\

(A21)~~~ $ F ( f ) ~=~ f $

\hfill $\Box$ \\

In view of (A17), it is obvious that condition (A21) is equivalent with \\

(A21$^*$)~~~ $ F ( f ) ( x ) \subseteq f ( x ),~~ x \in \Omega $ \\

Furthermore, (A20), (A21) give \\

(A22) \quad $ F ( f ) ~~\mbox{is~ s-continuous~ for}~~ f \in {\bf A} ( \Omega ) $ \\ \\

{\bf Example A1} \\

Let us illustrate the concept of s-continuity in the simplest case of functions with one
variable and one single discontinuity. Thus, with $\Omega = {\bf R}$, we take $f : \Omega
\longrightarrow \overline{{\bf R}}$, or in other words, $f \in {\cal A} ( \Omega )$, defined
by \\

$ f ( x ) ~=~ \begin{array}{|l}    ~~a ~~~\mbox{if}~~ x < 0 \\ \\
                                                 ~~b ~~~\mbox{if}~~ x = 0 \\ \\
                                                 ~~c ~~~\mbox{if}~~ x > 0
                    \end{array} $ \\ \\

where $a, b, c \in \overline{{\bf R}},~ a \neq c$. Then $f$ is {\it not} s-continuous. \\

Let us now take $f : \Omega \longrightarrow \overline{{\bf I R}}$, that is,  $f \in {\bf A}
( \Omega )$, defined by \\

$ f ( x ) ~=~ \begin{array}{|l}    ~~a ~~~\mbox{if}~~ x < 0 \\ \\
                                                 ~~[~ b,~ c ~] ~~~\mbox{if}~~ x = 0 \\ \\
                                                 ~~d ~~~\mbox{if}~~ x > 0
                    \end{array} $ \\ \\

where $a, b, c, d \in \overline{{\bf R}},~ a \leq d$ and $b \leq c$. Then f is s-continuous,
if and only if $b \leq a$ and $d \leq c$. Similarly, if $a \geq d$ and $b \leq c$, then f is
s-continuous, if and only if $b \leq d$ and $a \leq c$. \\

Consequently, returning to the first example above, it follows that f is s-continuous, if and
only if $a = b = c$, that is, if and only if f is continuous, see (A2.1.32) below, for the
general case of functions $f \in {\cal A} ( \Omega )$.

\hfill $\Box$ \\

The {\it fundamental} concept is presented now in \\ \\

{\bf Definition A2} \\

A function $f \in {\bf A} ( \Omega )$ is called {\it Hausdorff-continuous}, or in short, {\it
H-continuous}, if and only if $f$ is s-continuous, and in addition, for every s-continuous
function $g \in {\bf A} ( \Omega )$, we have satisfied the {\it minimality} condition on
$f$ : \\

(A23)~~~ $ g ( x ) \subseteq f ( x ),~~ x \in \Omega ~~~\Longrightarrow~~~  g ~=~ f $ \\

We shall denote by ${\bf H} ( \Omega )$ the set of all Hausdorff-continuous interval valued
functions on $\Omega$. \\ \\

{\bf Example A2} \\

Let us again consider the second function in Example A1 above. Then $f$ is H-continuous, if
and only if $a = b$ and $c = d$ \\

Let us give three further examples. \\

First, let us define $\alpha : {\bf R} \longrightarrow \overline{{\bf IR}}$ by \\

$ \alpha ( x ) ~=~ \begin{array}{|l} ~ - 1 ~~~~\mbox{if}~~ x < 0 \\ \\
                                                      ~~~ [~ - 1,~ 1 ~] ~~~~\mbox{if}~~ x = 0 \\ \\
                                                      ~~~ 1 ~~~~\mbox{if}~~ x > 0
                            \end{array} $ \\ \\

and then we can define $\beta : {\bf R}^2 \longrightarrow \overline{{\bf IR}}$ by \\

$ \beta ( x,~ y ) ~=~ \begin{array}{|l} ~~ \alpha ( \sin ( 1 / ( x^2 + y^2 ) ) ) ~~~~\mbox{if}~~ ( x,~ y )  \neq ( 0,~ 0 ) \\ \\
                                                          ~~~ [~ - 1,~ 1 ~] ~~~~\mbox{if}~~ ( x,~ y ) = ( 0,~ 0 )
                                \end{array} $ \\ \\

It is easy to see that both $\alpha$ and $\beta$ are H-continuous. \\

The third example is a typical {\it shock wave} solution of the well known nonlinear PDE in
Fluid Dynamics \\

$ U_t + U U_x ~=~ 0,~~~ t \geq 0,~ x \in {\bf R} $ \\

which corresponds to the initial value problem \\

$ U ( 0,~ x ) ~=~ \begin{array}{|l} ~~ 1 ~~~~\mbox{if}~~ x \leq - 1 \\ \\
                                                     ~ - x~~~~\mbox{if}~~ - 1 \leq x \leq 0 \\ \\
                                                     ~~ 0 ~~~~\mbox{if}~~ x \geq 0
                           \end{array} $ \\ \\

Namely, with $\Omega = [ 0, \infty ) \times {\bf R}$, we have the solution $U : \Omega
\longrightarrow \overline{{\bf IR}}$ given by \\

$ U ( t,~ x ) ~=~ \begin{array}{|l} ~~ 1 ~~~~~~\mbox{if}~~~~ 0 \leq t < 1,~~ x < t - 1 \\ \\
                                                    ~~ x / ( t - 1 ) ~~~~~~\mbox{if}~~~~ 0 \leq t < 1,~~ t - 1 \leq x \leq 0 \\ \\
                                                    ~~ 0 ~~~~~~\mbox{if}~~~~ 0 \leq t < 1,~~ x > 0 \\ \\
                                                    ~~ 1 ~~~~~~\mbox{if}~~~~ t \geq 1,~~ x < ( t - 1 ) / 2 \\ \\
                                                    ~~ [~ - 1,~ 1 ~] ~~~~~~\mbox{if}~~~~ t \geq 1,~~ x = ( t - 1 ) / 2 \\ \\
                                                    ~~ 0 ~~~~~~\mbox{if}~~~~ t \geq 1,~~ x > ( t - 1 ) / 2
                           \end{array} $ \\ \\

Then $U$ is H-continuous. \\ \\

{\bf Remark A1} \\

The {\it minimality} condition (A23) in the
above definition of H-continuous functions proves to play a fundamental role.

\hfill $\Box$ \\

As for the significance of the {\it regularity} property of being Hausdorff continuous, here
we an important {\it similarity} between usual continuous, and on the other hand,
Hausdorff-continuous functions, on the other. Namely, both of them are determined {\it
uniquely} if they are known on a {\it dense} subset of their domains of definition. \\
This property comes in spite of the fact that Hausdorff-continuous functions can have
discontinuities on sets of first Baire category, and such sets can have arbitrary large
positive Lebesgue measure. \\

Indeed, we have \\

{\bf Theorem A1} \\

Let $f = [~ \underline{f},~ \overline{f} ~],~ g = [~ \underline{g},~ \overline{g} ~] \in
{\bf A} ( \Omega )$ be two H-continuous functions, and suppose given any dense subset
$D \subseteq \Omega$. Then with the partial order in (A5), we have \\

a)~ $ \underline{f} ( x ) ~\leq~ \underline{g} ( x ),~~  x \in D ~~~\Longrightarrow~~~ f ~\leq~ g ~~\mbox{on}~~ \Omega $ \\

b)~ $ \overline{f} ( x ) ~\leq~ \overline{g} ( x ),~~ x \in D ~~~\Longrightarrow~~~ f ~\leq~ g ~~\mbox{on}~~ \Omega $ \\

c)~  $ f ( x ) ~\leq~ g ( x ),~~ x \in D ~~~\Longrightarrow~~~ f ~\leq~ g ~~\mbox{on}~~ \Omega $ \\

Also \\

d)~ $ f ( x ) ~=~ g ( x ),~~ x \in D ~~~\Longrightarrow~~ f ~=~ g ~~\mbox{on}~~ \Omega $ \\ \\

{\large \bf References} \\

1. Anguelov R, Rosinger E E : Solution of nonlinear PDEs by Hausdorff Continuous Functions (to
appear) \\

2. Arnold V I : Lectures on PDEs. Springer Universitext, 2004 \\

3. Grosser M, et.al. : Geometric Theory of Generalized Functions with Applications to General
Relativity. Mathematics and its Applications, Vol. 573, Kluwer, Dordrecht, 2001 \\

4. Oberguggenberger M B, Rosinger E E : Solutions of Continuous Nonlinear PDEs through Order
Completion. North-Holland Mathematics Studies, Vol. 181, (432 pages). Amsterdam, 1994 \\

5. Lewy H : An example of a smooth linear partial differential equation without solution. Ann.
Math., Vol. 66, No. 2, 1957, 155-158 \\

6. Rosinger E E : Parametric Lie Group Actions on Global Generalized Solutions of Nonlinear
PDEs, including a Solution to Hilbert's Fifth Problem, (234 pages). Kluwer, Dordrecht, 1998 \\

7. Rosinger E E : How to solve smooth nonlinear PDEs in algebras of generalized functions with
dense singularities (invited paper). Applicable Analysis, Vol. 78, 2001, 355-378 \\

\end{document}